# Max/Min Puzzles in Geometry†

Revised and Corrected

James M Parks
Dept. of Math.
SUNY Potsdam
Potsdam, NY
*parksjm@potsdam.edu*
January 4, 2022

**Abstract.** The objective here is to find the maximum polygon, in area, which can be enclosed in a given triangle, for the polygons: parallelograms, rectangles and squares. It is initially assumed that the choices are *inscribed polygons*, that is *all vertices of the polygon are on the sides of the triangle*. This concept is generalized later to include *wedged polygons*.

One of the earliest examples of a max/min inscribed polygon puzzle is known as *Euclid's Maximum Problem* [7], stated here as Puzzle 1.

**Puzzle 1.** *Inscribe a parallelogram in a given triangle, with its sides parallel to two sides of the triangle, such that the parallelogram has maximum area.*

First question: *"How do we solve this Puzzle?"*. Since all vertices of the parallelogram must be on the triangle, and the sides must be parallel to *2* sides of the triangle, we must have a situation similar to one of the parallelograms *BDEF* or *CHIG* or other in $\Delta ABC$, Fig. 1.

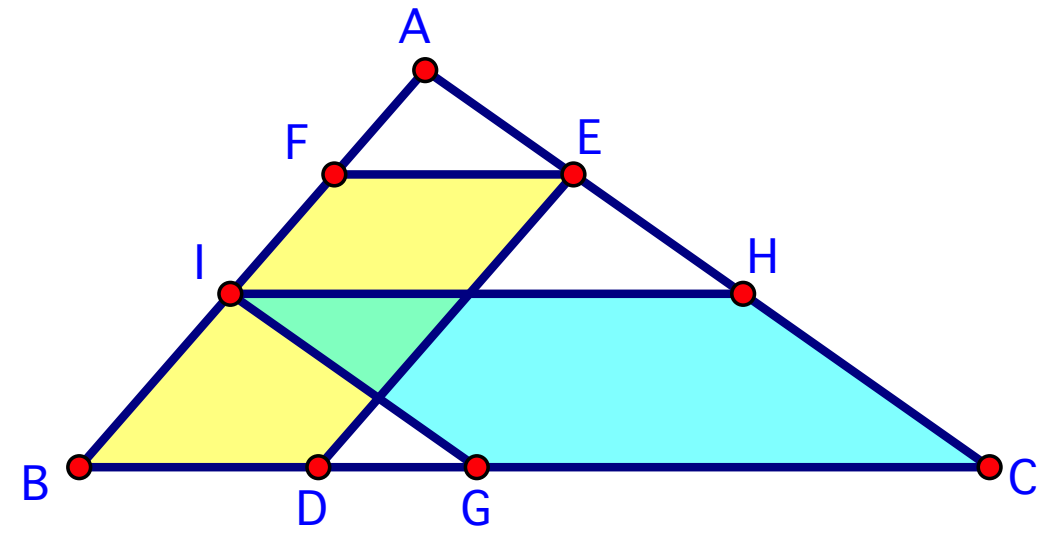

*Figure 1*

Using *Sketchpad*,* or similar, makes the study of such examples a simple process. It then becomes apparent that one of the prime candidates for the maximum parallelogram for the configuration of *BDEF* looks like Fig. 2, where the vertices *D, E, F* are very near the midpoints of sides *BC, CA,* and *AB,* respectively.

---

† Some of these results have appeared previously [4] , [5].

* All graphs in this paper were made with *Sketchpad v.5.10 BETA*.

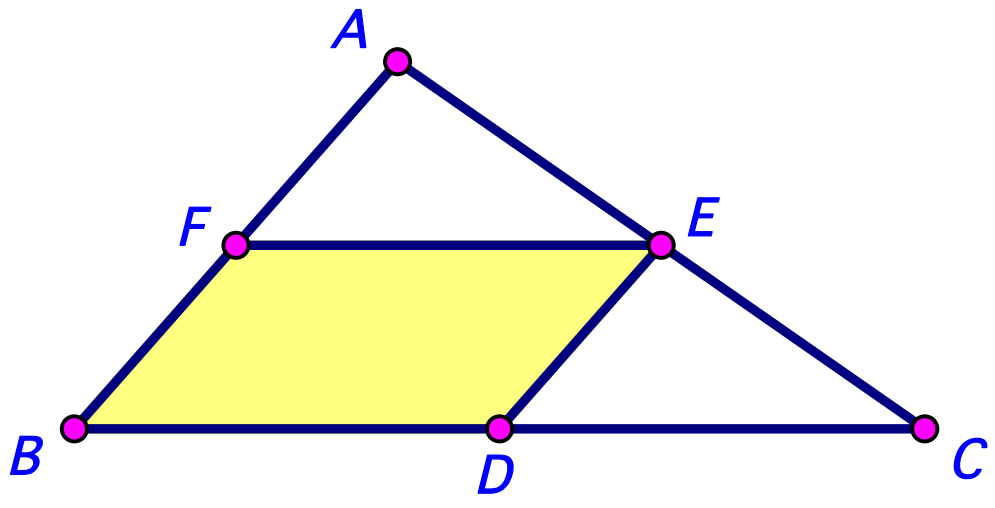

*Figure 2*

It can then be conjecture that the maximum inscribed parallelograms are determined by the midpoints of the sides of the triangle and the vertices of the triangle. A proof of this conjecture follows from Euclid's Proposition *27* in Book VI of *The Elements* [3], but there exists also a proof which is based on algebra, which we give here.

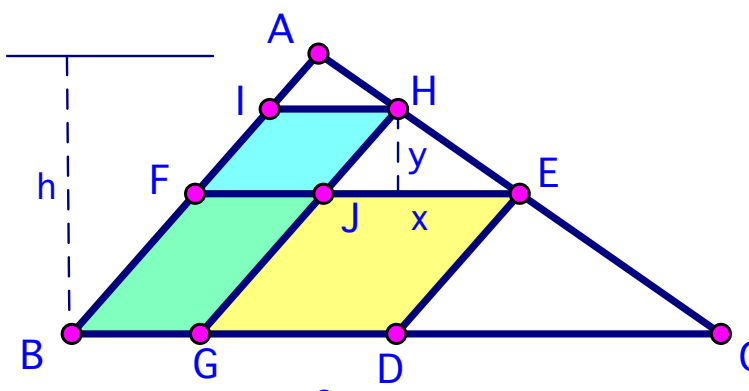

*Figure 3*

Let $h$ be the height of $\Delta ABC$ at $A$, and $a$ the length of the base $BC$, Fig.3.
Then $(ABC)=ha/2$, (parentheses denote area).
Assume $BDEF$ is the desired parallelogram with $D$, $E$, $F$ the midpoints of the sides of $\Delta ABC$ as shown.
Then $(BDEF) = (h/2)(a/2) = ha/4$.
Let $BGHI$ be another inscribed parallelogram with side $BI$ on $AB$, and let $x$ equal the length of $JE$, $J$ the intersection of $GH$ with $FE$, and $y$ the height of $\Delta JEH$ at $H$, Fig. 3.
Then $(BGHI) = (a/2 - x)(h/2 + y) = ha/4 - hx/2 + ay/2 - xy$.
But $\Delta ABC \sim \Delta HJE$, since the angles are equal, so $x/y = a/h$, and thus $hx = ay$.
Therefore, $(BGHI) = (BDEF) - xy$, and $BDEF$ is clearly the maximum parallelogram.
Similarly, if one chooses $H$ on $EC$, and $G$ on $DC$, and similarly if one chooses vertex $C$ or vertex $A$. Thus the max parallelogram occurs when *3* of the vertices are midpoints of the triangle sides, and the forth vertex is a vertex of the triangle.
Hence there exist *3* max parallelogram solutions for $\Delta ABC$ regardless of whether $\Delta ABC$ is acute, right, or obtuse.
There are obviously no min parallelograms.
These *3* parallelograms are all equal in area, since any two of them share a base and a height, as can be clearly seen in Fig. 4.

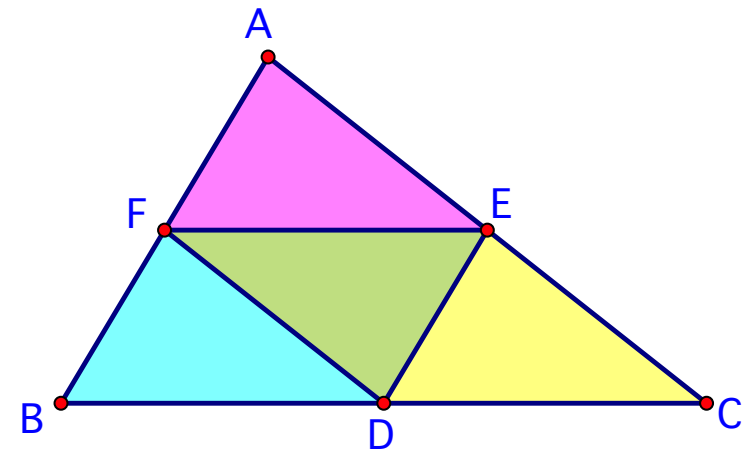

*Figure 4*

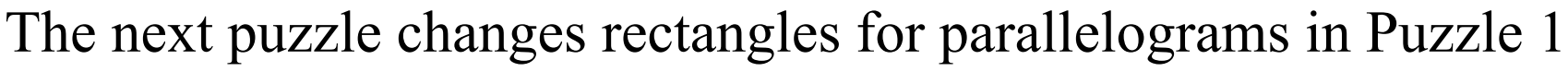

The next puzzle changes rectangles for parallelograms in Puzzle 1.

**Puzzle 2.** *Inscribe a rectangle in a given triangle, such that the rectangle has maximum area.*

To begin with the relationship between inscribed rectangles and inscribed parallelograms on the same vertices *E, F* of sides *AB, AC,* resp., is illustrated in Fig. 5, where the base angles are acute. It is well known that inscribed maximum area rectangles must have one side on a side of the triangle [2].

Also, the maximum rectangle must 'coincide' with the maximum parallelogram, in the sense shown in Fig. 5, since if not, then there would be a larger rectangle or parallelogram, which would contradict the maximum parallelogram or rectangle, as the case may be. Thus the maximum rectangle area is achieved when the side opposite the side which is on $\Delta ABC$ is connecting the midpoints of the other *2* sides of $\Delta ABC$.

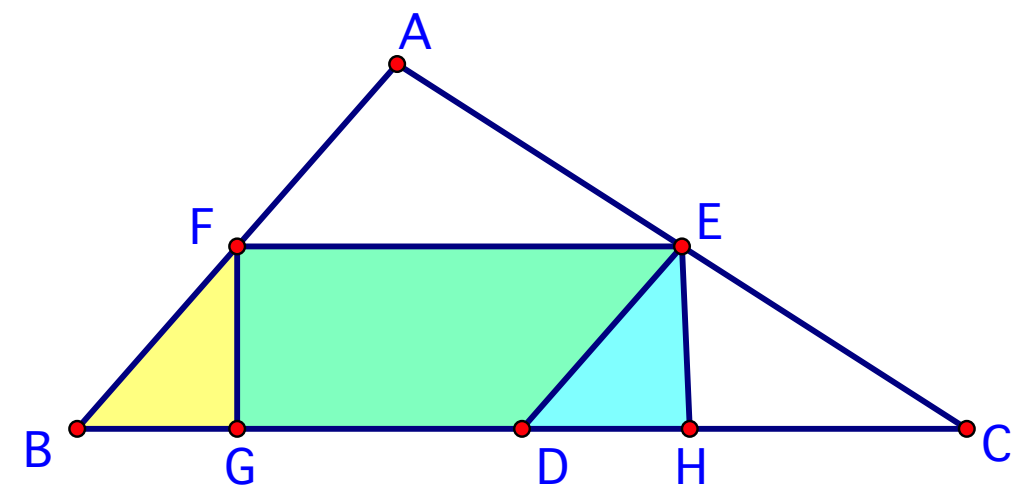

*Figure 5*

**Example 1.** For example, let $\Delta ABC$ be the triangle in the *xy*-plane, where $A = (5,3)$, $B = (0,0)$, and $C = (6,0)$, Fig.6. Then $(ABC) = 3(6/2)u^2 = 9u^2$, since the height $h = 3$, and the base $a = 6$. Let *D* be the midpoint of *AC, E* the midpoint of *AB*, and *F* the midpoint of *BC*. So $D = (5.5,1.5)$, $E = (2.5,1.5)$, and $F = (3,0)$. Then the parallelogram has area $(CDEF) = (h/2)(a/2) = (3/2)3u^2 = 9/2u^2$. Also, if *DH* is perpendicular to *BC* at $H = (1.5,0)$, and *EG* is perpendicular to *BC* at $G = (2.5,0)$, then the rectangle has area $(HDEG) = (CDEF) = 9/2u^2$.

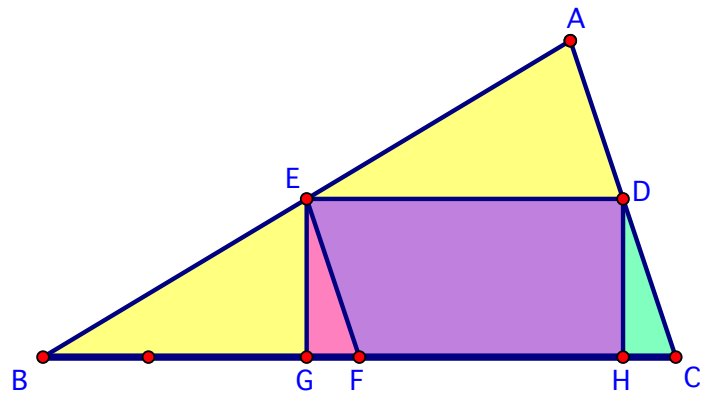

*Figure 6*

However, there is one difference between these two types of polygons. For the parallelograms, there are *3* equal maximum ones in all types of triangles. But for rectangles, there are *3* equal maximum rectangles in acute triangles, *2* equal maximum rectangles in right triangles (the max rectangles on sides *AB* and *AC* coincide), and one maximum rectangle in obtuse triangles, Fig. 7. There are obviously no min rectangles. Also, you can't inscribe a rectangle in a triangle if one of the base angles is obtuse.

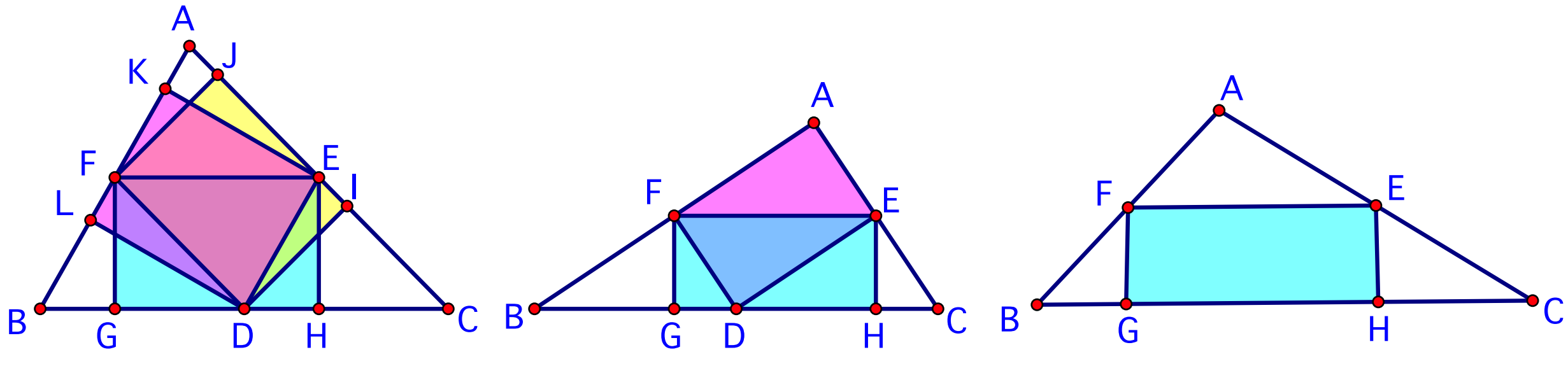

*Figure 7*

As with parallelograms, the area of each of the inscribed rectangles is *ha/4,* where *h* is the height, and *a* is the length of the respective base side. So in the acute triangle in Fig. 7, *(GHEF)* = *(DIJF)* = *(KLDE)* = *ha/4,* and similarly for right triangles, and obtuse triangles.

Since we are studying inscribed rectangles in triangles, we can also consider the case of inscribed squares. However, this turns out not to be just a special case of rectangles.
But the first problem is "how do you inscribe a square in a triangle"?
With a rectangle the height and the width are independent, not so with a square.

Here's a method, due to Polya [6], for constructing an inscribed square on side *BC* of a given triangle Δ*ABC* with no obtuse angle on the base. Of course there are also other methods [2].

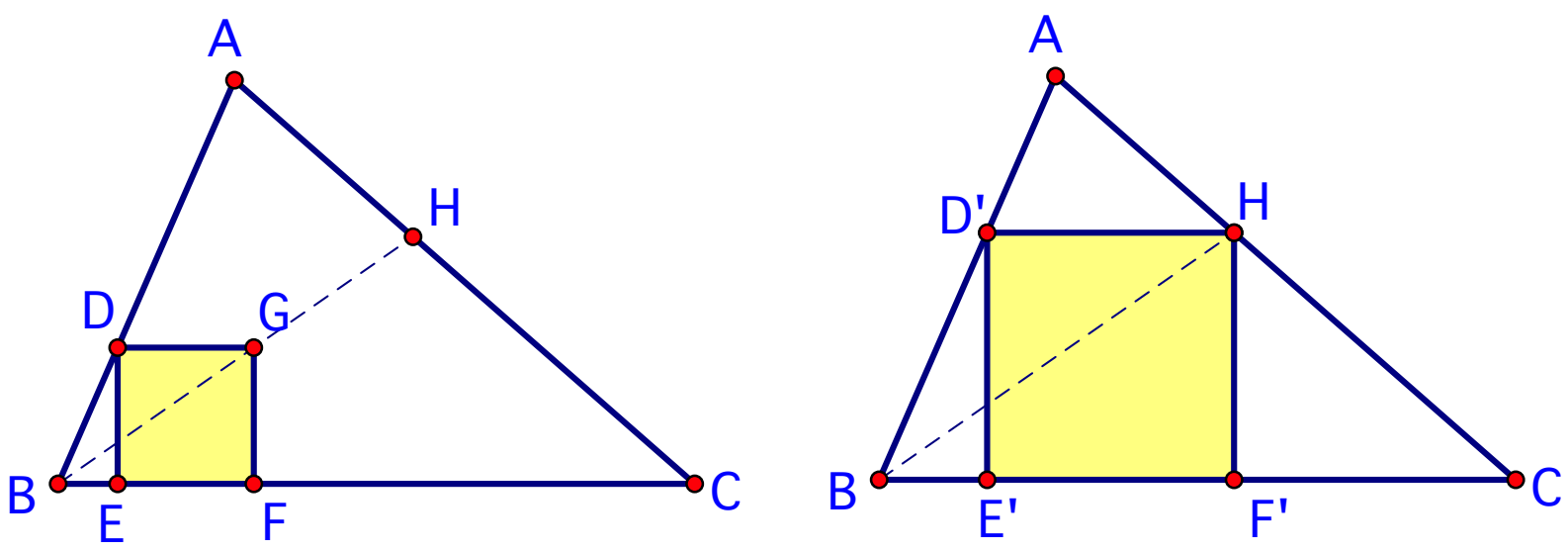

*Figure 8*

Let *D* be a point on side *AB*, near *B,* and construct a square *DEFG* with *EF* on side *BC,* Fig. 8 Left. Construct a line from *B* through *G* to *AC* at *H.* The desired inscribed square *D'E'F'H* is obtained by the dilation of *DEFG* about point *B* by the ratio *BH/BG,* Fig. 8 Right.

As with rectangles, the number of inscribed squares depends on the type of triangle. There are *3, 2,* or *1* inscribed squares, if the triangle is acute, right, or obtuse, respectively.
However, unlike rectangles, the area of the squares may differ with the base side of the triangle, since the height equals the width, Fig. 9.

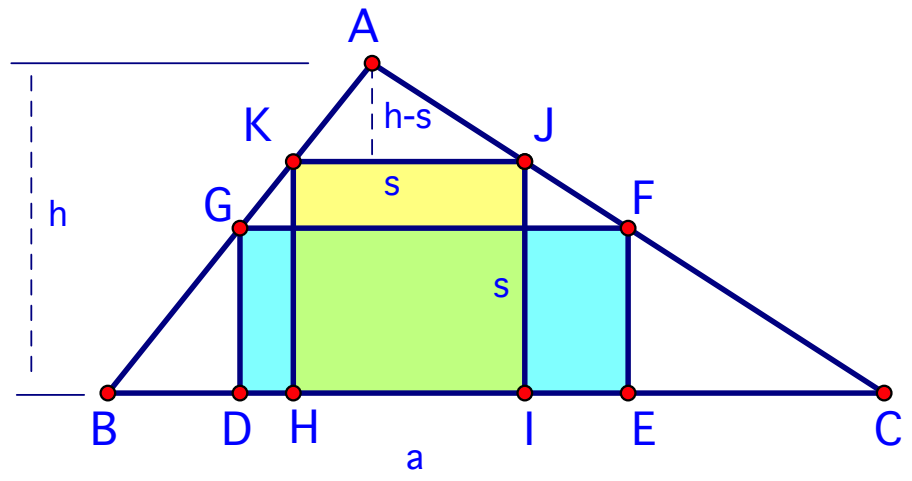

*Figure 9*

**Puzzle 3.** *Determine a formula for the length of the side of an inscribed square.*

For the formula for computing the length *s* of the side of an inscribed square proceed as follows [2]. Let Δ*ABC* be a given triangle with height *h,* and inscribed square *HIJK* on base *BC* of length *a,* where the base angles are acute, Fig. 9. Then, let *s* be the length of the side of the square, and note that Δ*ABC* ~ Δ*AKJ,* (equal angles). This determines the equation, s/a = (h-s)/h.
Solving for *s* determines the solution *s = ha/(h+a).*
Thus, the area of the inscribed square with side *s* is $s^2 = [ha/(h+a)]^2$.

Notice that the area of the inscribed square seems to be smaller than the area of the maximum inscribed rectangle on the same side, Fig. 9. This leads us to the new Puzzle 4.

**Puzzle 4.** *Show that the area of a maximum inscribed rectangle is always larger than the area of the inscribed square on the same side of a given triangle which has no obtuse angles on the base, with one exception, if the height equals the base the areas are equal.*

Consider $\Delta ABC$ with inscribed square *HIJK,* and the inscribed rectangle *DEFG*, Fig. 9.
Then we need to show that $(ha/(h+a))^2 \leq ha/4$.
But this inequality is equivalent to the inequality $0 \leq (h - a)^2$, which is obviously true, with the added bonus that equality holds only when $h = a$.
Hence the puzzle is solved.

**Example 2.** When looking at examples of acute triangles notice that the largest inscribed square is always on the smallest side, and the smallest inscribed square is always on the largest side.
Let $\Delta ABC$ be the acute triangle with $m<A = 75°$, $m<B = 60°$, and $m<C = 45°$, Fig. 10. Then $a > b > c$, since the largest side is opposite the largest angle.
Let $c = 2$, then altitude $h_a = \sqrt{3} \sim 1.73$, $a = 1+\sqrt{3} \sim 2.73$, and $s_a = \sqrt{3}(1+\sqrt{3})/(1+2\sqrt{3}) \sim 1.060$; $b = \sqrt{6} \sim 2.45$, $h_b \sim 1.93$, so $s_b \sim 1.080$; and $h_c \sim 2.37$, so $s_c \sim 1.084$.
Thus $s_a < s_b < s_c$.

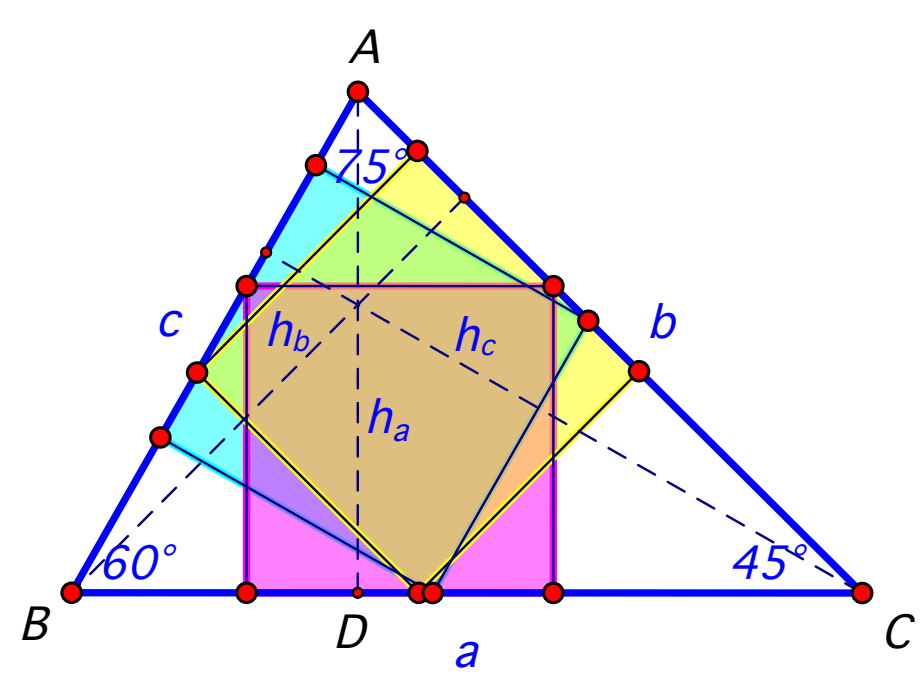

*Figure 10*

This observation is the basis for Puzzle 5.

**Puzzle 5.** *If $\Delta ABC$ is an acute scalene triangle, with sides $a > b > c$, show that the sides of the inscribed squares satisfy $s_a < s_b < s_c$.*

It is only necessary to show that whenever the triangle sides satisfy $a > b$, then the square sides satisfy $s_a < s_b$.
So, let $\Delta ABC$ is an acute triangle, with sides which satisfy $a > b$, and $h_b$ is the height of $\Delta ABC$ at *B,* Fig. 11.
First we show that $a + h_a > b + h_b$.
By the assumption, $a - b > 0$, and $h_a < h_b$, since $(ABC) = ah_a/2 = bh_b/2$. Plus it follows that $b > h_a$, since *b* is the hypotenuse of a right triangle with one leg of length $h_a$.
Thus $b - h_a > 0$, so $(a + h_a) - (b + h_b) = (a - b) + (2(ABC)/a - 2(ABC)/b) = (a - b)(1 - 2(ABC)/ab) = (a - b)(b - h_a)/b > 0$.
Then $s_a = ah_a/(a+h_a) < ah_a/(b+h_b) = bh_b/(b+h_b) = s_b$, so $s_a < s_b$.

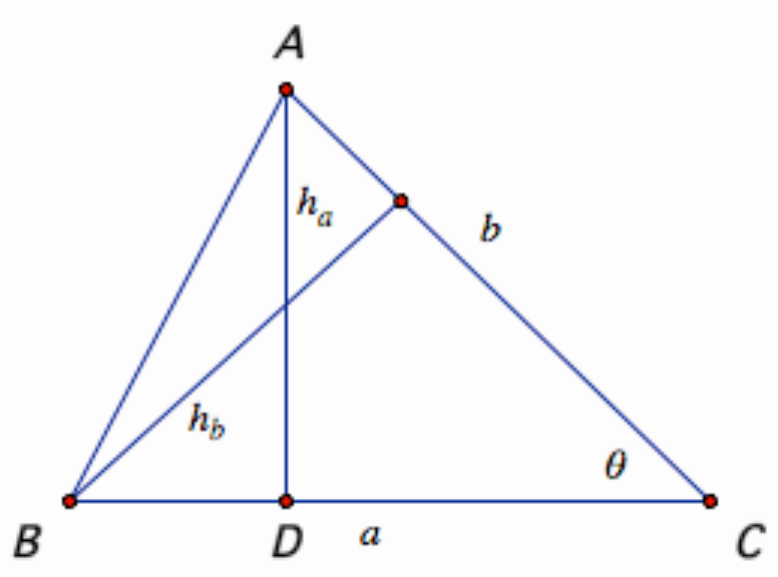

*Figure 11*

Similar arguments show that if $b > c$, then $s_b < s_c$.
Thus the maximum inscribed square is on the shortest side, and the minimum inscribed square is on the longest side.

Clearly an equilateral triangle has *3* equal inscribed squares.
But what about the obtuse triangles?
The problem with obtuse triangles is that if there is a square at the obtuse angle, there is no way that the square can be inscribed in the triangle, that is all of the squares vertices cannot be on the triangle, Fig.12. Thus *3* vertices of the enclosed square are in the triangle.
So instead of inscribed squares, we have something more general, called *wedged squares,* named by E. Calabi [1], [2].

There is a short-cut to constructing such squares at an obtuse angle, such as at vertex *A* in Fig.12.
Construct a line at the vertex *A* at a *45°* angle to the base leg *CA.* The intersection point *D* of this line with side *BC* determines a diagonal of the square *AFDE.*

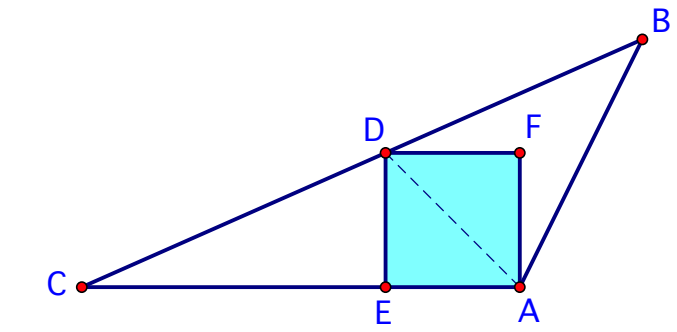

*Figure 12*

Of course we are not surprised that the area of wedged squares is computed by a different formula than those above. This determines Puzzle 6.

**Puzzle 6.** *Let ΔABC be an obtuse triangle, with obtuse angle at vertex A. Let AFDE be the largest wedged square on side CA = b. Find a formula for the length of the side of the square* $s_b$ *in terms of ø, the angle at C.*

By the given, if ø is the angle at *C,* and $s_b$ is the length of the side of the square on side *b,* then $s_b/(b - s_b) = tanø = sinø/cosø$.
Solve this for $s_b$, we find $s_b = bsinø/(sinø + cosø)$, Fig.12.

Using wedged squares, there will be *3* wedged squares in every type of triangle, as we found with parallelograms above. The areas of wedged squares differ with the length of the base side.

The concept of a wedged square can also be applied to rectangles in an obtuse triangle Δ*ABC,* Fig.13. Given a rectangle on a base side *AC* with an obtuse angle at *A*, let one vertex of the rectangle be on *A*, and another vertex *F* on side *AC*. Then, using Sketchpad, vary the vertex *E* on side *BC.* The max area of the rectangle *(ADEF)* appears to be when vertex *F* is near the midpoint of side *AC*. We conjecture that this is the actual max area of a wedged rectangle in an obtuse triangle on the side with an obtuse angle, see Puzzle 7.

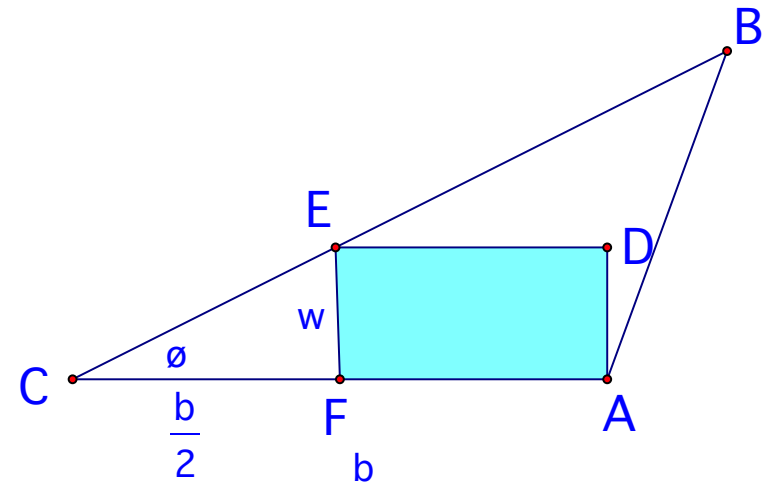

*Figure 13*

**Puzzle 7.** *If* $\Delta ABC$ *is an obtuse triangle, with obtuse angle at A, show that the area of a maximum wedged rectangle ADEF on side AC, with one vertex at A, and one vertex E on c = BC, occurs when the length of the base AF = b/2, for b = AC.*

Assume *ADEF* is the wedged rectangle on *AC* with *F* on the midpoint of *AC*. Then $(ADEF) = (b^2/4)Tan\text{ø}$, ø = angle *C,* Fig. 13.
Suppose *AD'E'F'* is another wedged rectangle on *AC,* with length $b/2 + e$, $e > 0$, and height $w' = (b/2 - e)Tan\text{ø}$, Fig. 14.
Solving the equation $(ADEF) \leq (AD'E'F')$ for *e,* determines $e = 0$.
Thus *ADEF* is the max wedged rectangle on side *AC.*
Similarly, for *AD'E'F'* with length $b/2 - e$, $e > 0$.

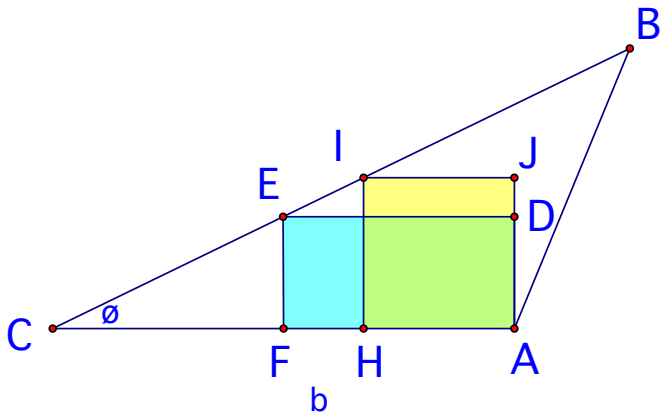

*Figure 14*

The areas of wedged rectangles depend on the length of the base side and the angle opposite the obtuse angle on the base side, so they vary with the side of the triangle, as seen in the case for wedged squares in Example 1 above.
However, it appears that the larger side has the larger area for wedged rectangles, contrary to the case for squares, Fig. 15.

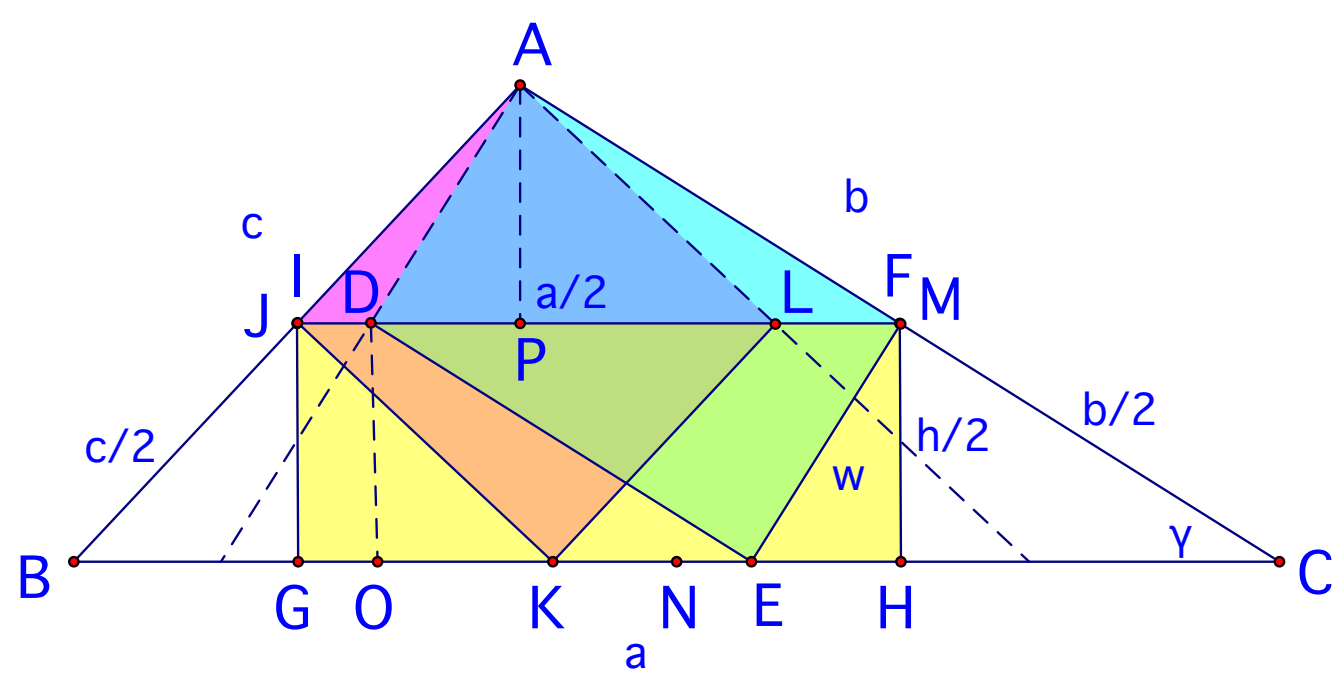

*Figure 15*

**Puzzle 8.** *Show that for the obtuse* $\Delta ABC$*, if* $a > b > c$*, with obtuse angle at A, then* $(GHFJ) \geq (ADEF) \geq (AJKL)$*, if F and J are on the midpoints M and I of AC and AB, respectively, for GHFJ the max rectangle on BC, and ADEF, AJKL wedged rectangles on AC, AB, respectively.*

That $(ADEF) \geq (AJKL)$ follows from the observation that angle *JAD* = angle *LAF*, so $FL > JD$, since $b > c$, Fig.15. Then, this is equivalent to $(ADF) > (AJL)$, since $(ADF) = (a/2 - JD)h/4$, and $(AJL) = (a/2 - LF)h/4$. But $(AJKL) = 2(AJL)$, and $(ADEF) = 2(ADF)$, so $(ADEF) \geq (AJKL)$.
To show $(GHFJ) \geq (ADEF)$, consider the triangle $\Delta ADF$, and the rectangle *GHFJ,* Fig.15. The rectangle *ADEF* is divided into *2* copies of $\Delta ADF$ by diagonal *DF,* one which is in *GHFJ,* since *D* is on *JF.* The vertical line *AP* cuts $\Delta ADF$ into *2* right triangles. One is $\Delta ADP$, which is congruent to $\Delta FEH$, with heights *h/2,* and sides *AD//FE,* since they are opposite sides of a rectangle, and *AP//FH,* since both are perpendicular to *DF*. The other one is $\Delta APF$ which is congruent to $\Delta DOE$, where $O$ is the foot of the vertical line at $D$ to $GH$. Again, the heights are $h/2$, and sides $AP//DO$, and $AF//DE$.
Thus, rectangle $DOHF$ equals $2$ copies of $\Delta ADF$, and $(DOHF) = (ADEF)$.
Since $GHFJ$ contains $DOHF$, it follows that $(GHFJ) > (DOHF) = (ADEF)$.

The relation between the areas of maximum wedged rectangles and wedged squares is the same as for rectangles and squares, Puzzle 4, Fig. 15.

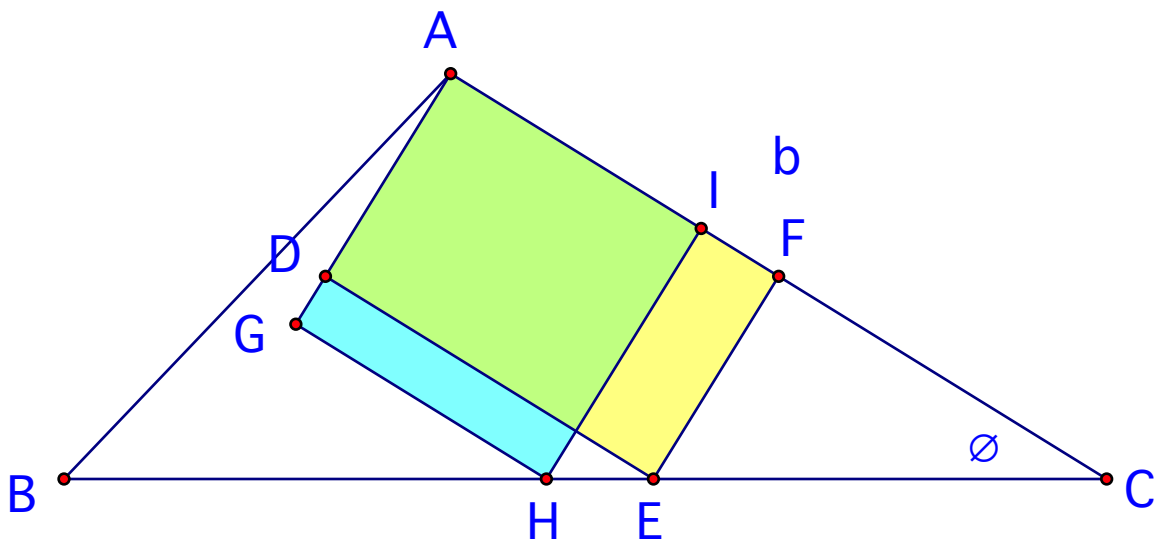

*Figure 16*

**Puzzle 9.** *Let $\Delta ABC$ be obtuse at A, and let ADEF be the maximum wedged rectangle on side AC with vertex F at the midpoint of AC, and AGHI the wedged square on AC. Show that $(ADEF) \geq (AGHI)$.*

Observe that the inequality $(ADEF) \geq (AGHI)$ is equivalent to $(b^2/4)Tan\text{ø} \geq (b\sin\text{ø}/(\sin\text{ø} + \cos\text{ø}))^2$, which is equivalent to $(\sin\text{ø} - \cos\text{ø})^2 \geq 0$, Fig.16. Equality holds at $\text{ø} = 45°$, because then $\Delta ABC$ is a right triangle with base $b$ equal to the height $c$.

It should be clear that the areas of wedged squares (and wedged rectangles) on the legs of isosceles triangles with central angle greater than $90°$ are equal. By allowing this more general type of enclosed squares in triangles, a new and unexpected result appears, Example 3. The best way to discover it is by using *Sketchpad.*

**Example 3.** Let $\Delta ABC$ be an obtuse isosceles triangle with obtuse angle at $A$, and assume sides $AB$ and $AC$ are fixed at length $2$. Construct the $3$ wedged squares on the sides of the triangle, then the area of the squares on the legs $AB$ and $AC$ are equal, Fig.17.

| Angle | Area of 2 eq. sq.s | Area of DEFG |
|---|---|---|
| 105.66° | 0.7438u$^2$ | 0.7677u$^2$ |
| 102.64° | 0.7908u$^2$ | 0.7968u$^2$ |
| 101.61° | 0.8070u$^2$ | 0.8061u$^2$ |
| 98.50° | 0.8569u$^2$ | 0.8325u$^2$ |
| 95.31° | 0.9095u$^2$ | 0.8565u$^2$ |

*Figure 17*

Now move the vertex *A* up and down in the range 95° between *105°*.
The chart in Fig.17 gives some values of the areas of the wedged squares as *A* is moved.
Observe that in the range *101.61°* to *102.64°* there is a reversal of the sizes from the *2* equal squares being larger than *DEFG*, to them being smaller. By continuity of the area function, there must be an angle, in the range *101.61°* to *102.64°*, where the *3* areas are all equal!
This result was discovered by E. Calabi, that is there exists a triangle, which is not equilateral, but which has *3* equal area wedged squares.
He calls it a *catastrophe,* in the sense of R. Thom [1].

It is possible to compute the exact values of the angle at *A*, and the ratio of the long side to the short side of $\Delta ABC$, Fig.17.

**Puzzle 10.** *Determine the angle at A, and the ratio of the long side to the short side of Calabi's triangle* $\Delta ABC$.

Assume $AC = AB = 1$. From above we have $s_a = ah_a/(a+h_a)$, and $s_c = s_b = b\sin ø/(\sin ø + \cos ø)$, where $ø$ = angle *C*. Then $\sin ø = h_a$ and $\cos ø = a/2$, so $s_a = s_b$ is equivalent to $a/(a+h_a) = 1/(h_a+a/2)$. Solving for $h_a$ determines that $h_a = (a^2 - 2a)/(2 - 2a)$.
Also, by the Pythagorean Theorem, $h_a = \sqrt{(1 - a^2/4)}$, so equating these two equations and squaring out the radical determines $(1 - a^2/4) = [(a^2 - 2a)/(2 - 2a)]^2$, which simplifies to the quartic equation $2a^4 - 6a^3 + a^2 + 8a - 4 = 0$. This equation has the extraneous root $a = 2u$, since *2* is a zero, but $a < 2$ by the assumptions. So, dividing out the factor $(a - 2)$, we obtain the cubic equation $2a^3 - 2a^2 - 3a + 2 = 0$.

The largest positive solution of this equation is approximately $a = 1.5513875...$, and since $b = 1$, this is the ratio of the longest side *a* to the shortest side *b* of $\Delta ABC$.
Since $\cos ø = a/2$, the angle $ø = 39.132...°$, approximately, and the obtuse angle at *A* is approximately *101.736...°*.

There are *8* different combinations of possible sizes of wedged squares in triangles, and they can be categorized as follows (cf. [8]).

Assume that the triangle $\Delta ABC$ satisfies the inequality $a \geq b \geq c$, with $a = 1$. Let $C = (0, 0)$, $B = (-1, 0)$, and $D = (-1/2, 0)$, Fig. 18.
There are *3* possible triangles shown here.
If angle *A* is a right angle, then it is on the circle with radius *a/2* and center *D.*
If angle *A* is obtuse, then *A* is inside this circle, and if angle *A* is acute, then *A* is outside this circle, but inside the circle centered at *C,* with radius *a*.

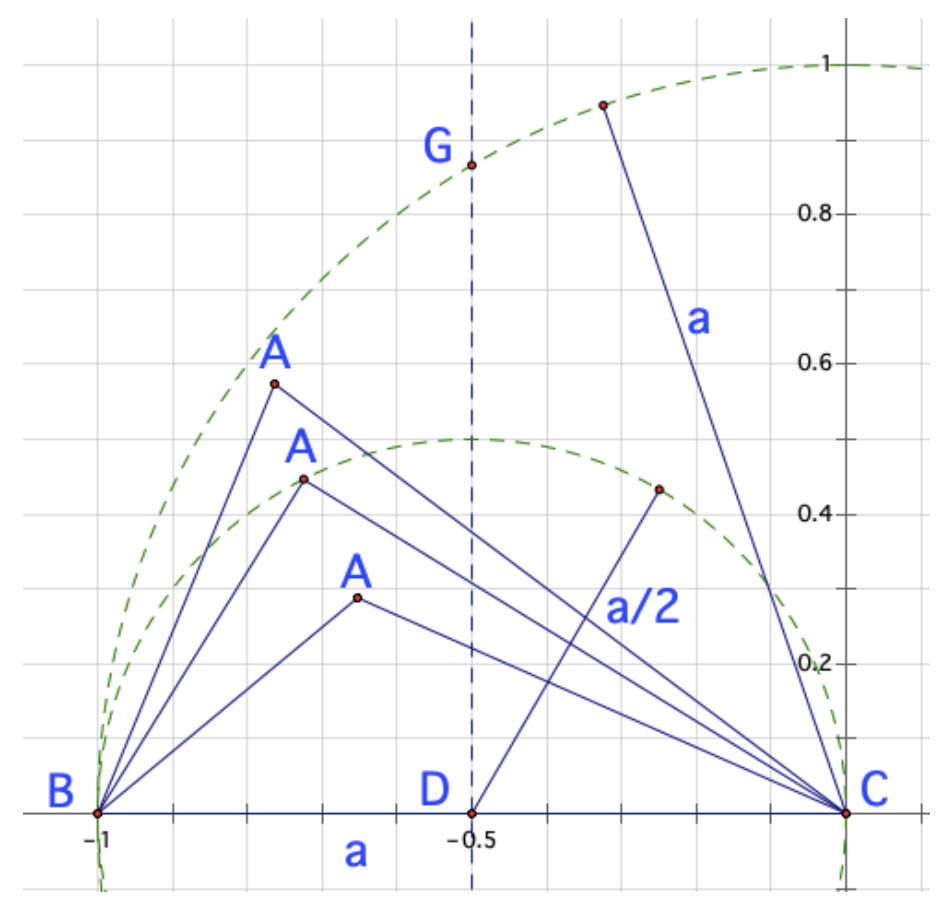

*Figure 18*

Whenever $A$ is on the vertical line $DG$, $\Delta ABC$ is isosceles, and if $A = G$, then $\Delta ABC$ is an equilateral triangle.
All possible points $A$ must be inside or on the circle about $C$ of radius $a$, since $a \geq b \geq c$.
If $s_a = s_b$, then $h_a = (\sin ø + \cos ø - 1)$, for $a = 1$, since $s_a = h_a/(1 + h_a)$, $s_b = b\sin ø/(\sin ø+\cos ø)$, and $h_a = b\sin ø$, $ø$ = angle $C$. Simplifying and formulating in polar coordinates $(r, \mu)$ determines the equation $r = (1 - \cot(\mu/2))$, where $\mu = 180° - ø$. Transforming this equation to rectangular coordinates determines the cubic equation $y^3 + x^2y + 2x^2 + 2y^2 + 2x = 0$, [8]. The graph of this equation contains the dashed curve $BHE$ under the semicircle centered at $D$, with radius $a/2$, and $s_b$ is the max for $A$ just above this curve, and $s_a$ is the max for $A$ just below it, Fig.19.
If $s_a = s_c$ we have the analogous case. This can be solved by reflecting $\Delta ABC$ about $DG$, solving for the cubic equation, then reflecting this curve about $DG$ again to get the solution for $\Delta ABC$.
The reflected cubic equation is $y^3 + (x+1)^2y + 2(x+1)^2 + 2y^2 - 2(x+1) = 0$. The graph of this equation contains the dashed curve $BIE$ between the semicircle and the cubic just below it. If $A$ is just above this curve, then $s_c$ is the larger, and $s_a$ is the minimum, and when $A$ is just below it, then $s_c$ is the minimum, and $s_a$ is the larger.
The two cubics intersect at point $E$, which is the point where Calabi's triangle occurs, Fig.19 [1].

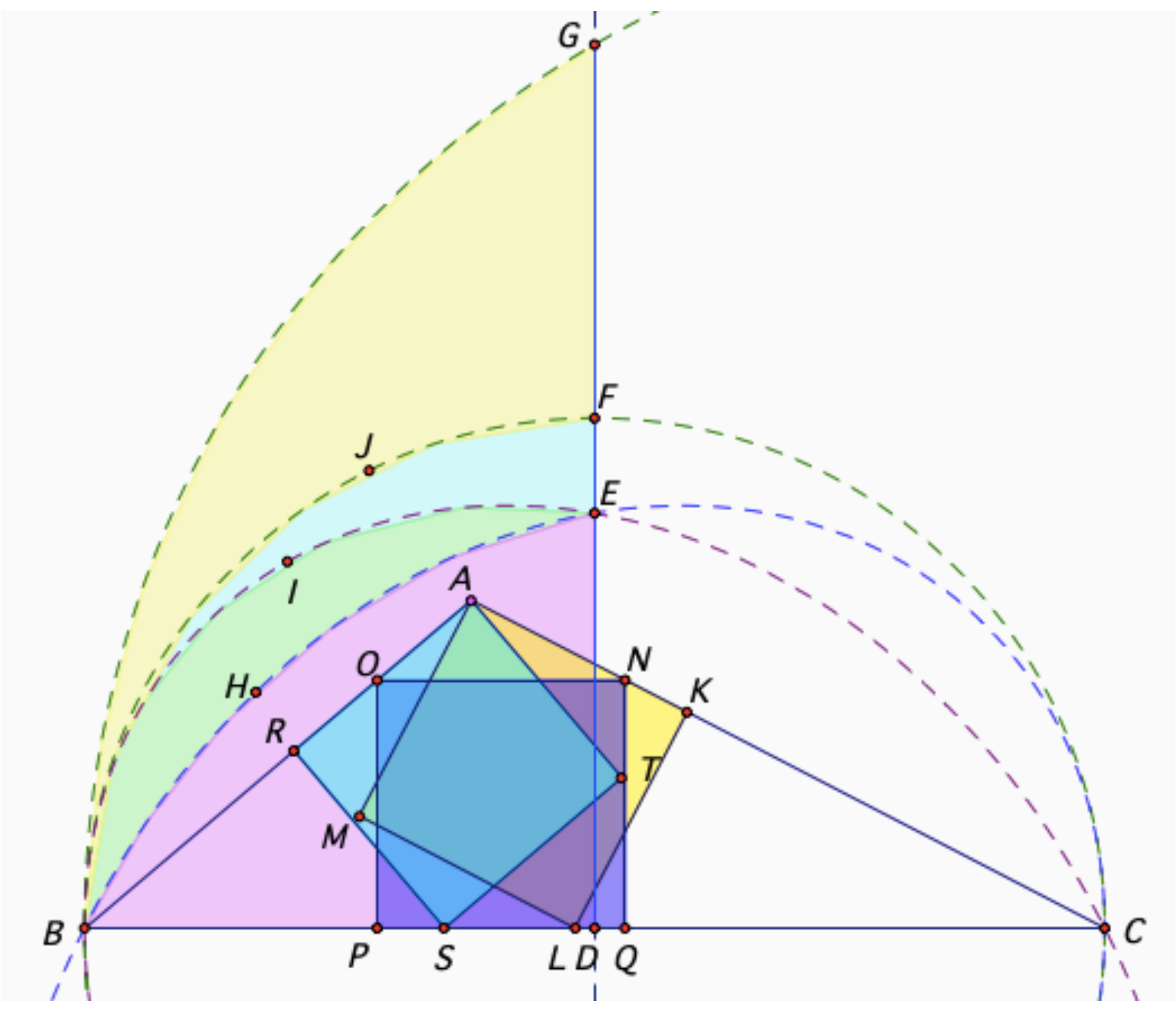

*Figure 19*

Starting at the top yellow region and moving down, Fig.19, the results are summarized in the following list.

$s_a < s_b < s_c$, for acute triangles: *A* in the yellow region,
$s_a < s_b = s_c$, for right triangles: *A* on circle *BJF,*
$s_a < s_c < s_b$, for obtuse triangles, and *A* in the blue region,
$s_a = s_c < s_b$, for *A* on the cubic curve *BIE,*
$s_c < s_a < s_b$, for obtuse triangles, and *A* in the green region,
$s_c < s_a = s_b$, for *A* on the cubic curve *BHE,*
$s_c < s_b < s_a$, for obtuse triangles, and *A* in the magenta region,
$s_a = s_b = s_c$, for Calabi's triangle, *A* on *E.*